\documentclass[11pt]{amsart}
\usepackage[T1]{fontenc}
\usepackage[utf8]{inputenc}
\usepackage{amstext}
\usepackage{amsthm}
\usepackage{amssymb}

\usepackage[hidelinks]{hyperref}

\makeatletter

\theoremstyle{plain}
\newtheorem{theorem}{\protect\theoremname}
\newtheorem{proposition}[theorem]{\protect\propositionname}
\newtheorem{lemma}[theorem]{\protect\lemmaname}
\newtheorem{remark}[theorem]{\protect\remarkname} 

\usepackage{babel}
\providecommand{\lemmaname}{Lemma}
\providecommand{\propositionname}{Proposition}
\providecommand{\theoremname}{Theorem}
\providecommand{\remarkname}{Remark} 

\begin{document}
\title[Rigidity of spacelike hypersurface]{Rigidity of spacelike hypersurface with constant curvature and intersection
angle condition}
\author{Shanze Gao}
\address{School of Mathematics\\ Guangxi University\\ Nanning 530004, Guangxi, People's Republic of China}
\address{Guangxi Center for Mathematical Research\\ Guangxi University\\ Nanning 530004, Guangxi, People's Republic of China}
\email{gaoshanze@gxu.edu.cn}


\begin{abstract}
In the Minkowski space, we consider a compact, spacelike hypersurface with boundary, which can be written as a graph on a spacelike hyperplane. We prove that, if its $k$-th mean curvature is constant, and its boundary is on the hyperplane with constant intersection angles, then the hypersurface must be a part of a hyperboloid, unless it is entirely contained in the hyperplane. The proof is based on an auxiliary function and associated integral equality.
\end{abstract}

\maketitle

\section{Introduction}

Characterization of hypersurfaces with constant curvature functions is a classical topic in differential geometry. In the Euclidean space, it has been known that the sphere is the only embedded, closed hypersurface with constant mean curvature or higher order mean curvature (see \cite{Aleksandrov,Montiel-Ros} etc.). Interested readers are referred to \cite{Hopf,Lopez} etc. for related results.

The Minkowski space $\mathbb{R}^{n,1}$ is $\mathbb{R}^{n+1}$ endowed with the Lorentzian metric \[\bar{g}=\sum_{i=1}^{n}dx_{i}^{2}-dx_{n+1}^{2}.\] Hypersurfaces with constant curvatures in $\mathbb{R}^{n,1}$ have attracted the interest of researchers. Cheng and Yau \cite{Cheng-Yau} prove the Bernstein type theorem for the maximal spacelike hypersurface. Treibergs \cite{Treibergs} shows that there are many entire spacelike hypersurfaces with constant mean curvature. Existence results for the cases of constant higher order mean curvatures are also obtained (see \cite{LiAnMin,Wang-Xiao} etc.). Among these hypersurfaces, the hyperboloid \[ x_{1}^2+\cdots+x_{n}^2-x_{n+1}^2=-1 \] is the umbilic one, which can be seen as the "sphere" in the Minkowski space.

It is shown that the hyperbolic caps and the hyperplanar balls are the only spacelike hypersurfaces with constant curvatures in \cite{Alias-Lopez-Pastor,Alias-Pastor98,Alias-Malacarne} etc., provided that the boundary of the hypersurface is spherical. Different with them, we consider an intersection angle boundary condition for the problem in this paper. 

Let $M$ be a connected, compact, spacelike hypersurface with boundary $\partial M$ and $\Pi$ be a spacelike hyperplane in $\mathbb{R}^{n,1}$. We define a function $\theta:\partial M\rightarrow \mathbb{R}$ by \[ \theta(p)=\bar{g}(N(p),N_{0})\quad\text{for }p\in\partial M,\] where $N$ and $N_{0}$ are the timelike unit normal vectors of $M$ and $\Pi$ respectively. We say that $M$ has \emph{constant intersection angles} if $\partial M\subset \Pi$ and $\theta$ is constant.

If there is a function $\tilde{u}$ defined on a domain $\overline{\Omega}\subset \Pi$ such that \[ M=\{x+\tilde{u}(x)N_{0}| x\in \overline{\Omega}\}, \] we say that $M$ is a graph on $\Pi$.

Let $\lambda=(\lambda_{1},...,\lambda_{n})$ denote the principal curvatures of $M$. The \emph{$k$-th mean curvature} (also called higher order mean curvature) of $M$ is
\[ H_{k}:=\frac{1}{\binom{n}{k}}\sigma_{k}(\lambda), \] where $\sigma_{k}(\lambda)$ is the $k$-th elementary symmetric polynomial of $\lambda$.

The following theorem shows that, the part of a hyperboloid cut by a spacelike hyperplane is the only spacelike hypersurface with constant $k$-th mean curvature and intersection angles, if $M$ is not entirely contained in $\Pi$.

\begin{theorem}\label{thm:main}
	Suppose $M$ is a compact, spacelike hypersurface with boundary $\partial M$ in $\mathbb{R}^{n,1}$, which is a graph on a spacelike hyperplane $\Pi$. If the $k$-th mean curvature of $M$ is constant, and $\partial M$ is on $\Pi$ with constant intersection angles, then $M$ is either a part of a hyperboloid, or entirely contained in $\Pi$.
\end{theorem}

The theorem can be seen as an overdetermined boundary value problem from a PDE perspective. In fact, the boundary conditions are equivalent to \[ \tilde{u}=0\text{ and }|D\tilde{u}|=c\quad\text{on }\partial\Omega. \] By using the moving plane method, Serrin \cite{Serrin} proved that the solutions to some elliptic equations are radially symmetric under such boundary conditions. Interested readers may see \cite{Weinberger,Wente,BNST,Jia,GMY} etc. for some related results.

To prove the theorem, we simplify the problem to the case  $\Pi=\mathbb{R}^{n}\times \{ c \}$ by an isometry of $\mathbb{R}^{n,1}$ (see Section \ref{sec:pre}). Then we show the hypersurface is $k$-convex in Section \ref{sec:convexity}. In the last section, we establish an integral equality and consider an auxiliary function via the maximum principle, which show that the hypersurface must be totally umbilical. The proof is inspired by Weinberger's approach \cite{Weinberger} to the Serrin's result.

\section{Preliminaries}

\subsection{Minkowski space}
\label{sec:pre}

Let $\langle\cdot,\cdot\rangle$ denote the scalar product in $\mathbb{R}^{n,1}$, i.e., \[ \langle Y,Z\rangle=y_{1}z_{1}+\cdots+y_{n}z_{n}-y_{n+1}z_{n+1} \] for $Y=(y_{1},...,y_{n},y_{n+1}),Z=(z_{1},...,,z_{n},z_{n+1})$.

Let $X:M^{n}\rightarrow \mathbb{R}^{n,1}$ be a connected, smooth hypersurface. It is called \emph{spacelike} if the induced metric $g$ of $M$ is Riemannian. A vector $Y\in \mathbb{R}^{n,1}$ is called \emph{timelike} if $\langle Y,Y \rangle<0$. It can be checked that the normal vector of a spacelike hypersurface must be timelike. A timelike vector $Y$ is called unit if $\langle Y,Y \rangle=-1$.

Let $\Pi$ be a spacelike hyperplane in $\mathbb{R}^{n,1}$. It can be written as \[ \Pi=\{ x\in \mathbb{R}^{n,1}|\langle x,N_{0}\rangle +c=0 \}, \] where $N_{0}$ is the timelike unit normal vector of $\Pi$ and $c$ is a constant. If a spacelike hypersurface $M$ is a graph on $\Pi$, then \[ M=\{x+\tilde{u}(x)N_{0}| x\in \overline{\Omega}\}, \] where $\tilde{u}$ is a smooth function on the closure of a smooth, bounded, connected open set $\Omega\subset \Pi$.

\begin{remark}\label{rem:u<0}
	If $M$ is not entirely contained in $\Pi$, we can assume that, there is a point $x_{0}\in \Omega$ such that $\tilde{u}(x_{0})<0$, by choosing a suitable $N_{0}$.
\end{remark}

\begin{remark}\label{rem:isomety}
	Let $O(n,1)$ be the group of all pseudo-orthogonal transformations of $\mathbb{R}^{n,1}$. Since $N_{0}$ is a fixed timelike unit vector, there exists a transformation $T\in O(n,1)$ such that $T(N_{0})=(0,...,0,1)$. Furthermore, $T(\Pi)=\mathbb{R}^{n}\times \{ c \}$ and $T(M)$ can be written as a graph on some domain in $\mathbb{R}^{n}\times \{ 0 \}$. Since $T$ is an isometry in $\mathbb{R}^{n,1}$, we know it maps a hyperboloid into another hyperboloid.
\end{remark}

Due to Remark \ref{rem:isomety}, we only need to prove Theorem \ref{thm:main} in the case $\Pi=\mathbb{R}^{n}\times \{ c \}$. From now on, we assume that the position vector $X$ of the hypersurface $M$ is in the form $X=(x,u(x))$, where $u(x)$ is a smooth function on the closure of a smooth, bounded, connected open set $\Omega\subset\mathbb{R}^{n}$. 

\subsection{Hypersurface in $\mathbb{R}^{n,1}$}
We recall some facts of hypersurfaces in the Minkowski space. Readers may refer to \cite{Cheng-Yau,LiAnMin} etc. 

Denote $E_{j}=(0,...,0,1,0,...,0)$, where $1$ is in the $j$-th place for $j\in\{1,2,...,n+1\}$. Then 
\[ \frac{\partial X}{\partial x_{i}}=E_{i}+u_{i}E_{n+1}. \]

The induced metric of $M$ is \[ g_{ij}=\langle\frac{\partial X}{\partial x_{i}},\frac{\partial X}{\partial x_{j}}\rangle=\delta_{ij}-u_{i}u_{j}. \] It is clear that $|Du|<1$ if $M$ is spacelike. Let $\nabla$ denote the Levi-Civita
connection on $M$. The Christoffel symbols are
\begin{equation*}
	\Gamma_{ij}^{m}=\frac{1}{2}g^{ml}(\frac{\partial g_{lj}}{\partial x_{i}}+\frac{\partial g_{il}}{\partial x_{j}}-\frac{\partial g_{ij}}{\partial x_{l}}),
\end{equation*}
where $(g^{ml})$ is the inverse matrix of $(g_{ij})$ and repeated indexes are summed from $1$ to $n$ here and later (unless otherwise stated). 

The timelike unit normal vector of $M$ is \[ N=\frac{(Du,1)}{\sqrt{1-|Du|^{2}}}. \]

Let $\nabla_{i}\nabla_{j}$ denote the Hessian operator. Then \[ \nabla_{i}\nabla_{j}X=\frac{\partial^{2}X}{\partial x_{i}\partial x_{j}}-\Gamma_{ij}^{m}\frac{\partial X}{\partial x_{m}}=h_{ij}N, \]
where $h_{ij}$ is the second fundamental form. 

It is easy to check that
\begin{equation*}
	h_{ij}=-\langle\frac{\partial^{2}X}{\partial x_{i}\partial x_{j}},N\rangle=\frac{u_{ij}}{\sqrt{1-|Du|^{2}}}.
\end{equation*}

The shape operator $A=(h_{i}^{j})$ satisfies
\begin{equation}\label{eq:dN}
	\frac{\partial}{\partial x_{i}}N=h_{i}^{j}\frac{\partial X}{\partial x_{j}}=h_{i}^{j}E_{j}+h_{i}^{j}u_{j}E_{n+1}.
\end{equation}

From \eqref{eq:dN}, we have
\begin{equation*}
	h_{i}^{j}=\frac{\partial}{\partial x_{i}}\langle N,E_{j}\rangle=\frac{\partial}{\partial x_{i}}\left(\frac{u_{j}}{\sqrt{1-|Du|^{2}}}\right)
\end{equation*}
and
\begin{equation}\label{eq:partialiN}
\frac{\partial}{\partial x_{i}}\langle N,E_{n+1}\rangle=-h_{i}^{j}u_{j}.
\end{equation}
It is clear that 
\[ \frac{\partial}{\partial x_{l}}h_{i}^{j}=\frac{\partial}{\partial x_{i}}h_{l}^{j}.\]
The covariant derivative of tensor $h_{i}^{j}$ satisfies
\[ \nabla_{j}h_{i}^{l}=\frac{\partial}{\partial x_{j}}h_{i}^{l}-\Gamma_{ij}^{m}h_{m}^{l}+\Gamma_{mj}^{l}h_{i}^{m}.\]
The Codazzi equation implies
\[ \nabla_{j}h_{i}^{l}=\nabla_{i}h_{j}^{l}. \]

\subsection{Symmetric functions of curvature}

Since principal curvatures $\lambda=(\lambda_{1},\dots,\lambda_{n})$
are the eigenvalues of the shape operator $A$, we have 
\[ \sigma_{k}(\lambda)=\sigma_{k}(A):=\frac{1}{k!}\delta_{j_{1}\cdots j_{k}}^{i_{1}\cdots i_{k}}h_{i_{1}}^{j_{1}}\cdots h_{i_{k}}^{j_{k}},\]
where $\delta_{j_{1}\cdots j_{k}}^{i_{1}\cdots i_{k}}$ is the generalized
Kronecker symbol defined by
\[
\delta_{j_{1}\cdots j_{k}}^{i_{1}\cdots i_{k}}:=\begin{cases}
	1, & \text{if }(i_{1}\cdots i_{k})\text{ is an even permutation of }(j_{1}\cdots j_{k}),\\
	-1, & \text{if }(i_{1}\cdots i_{k})\text{ is an odd permutation of }(j_{1}\cdots j_{k}),\\
	0, & \text{otherwise}.
\end{cases}
\]
We also set $\sigma_{0}(A)=1$ and $\sigma_{k}(A)=0$ for $k<0$ and
$k>n$. The $k$-th mean curvature of $M$ can be written by 
\[ H_{k}=\frac{1}{\binom{n}{k}}\sigma_{k}(A). \]

Define
\[(\sigma_{k}(A))_{j}^{i}:=\frac{\partial\sigma_{k}(A)}{\partial h_{i}^{j}}=\frac{1}{(k-1)!}\delta_{j_{1}\cdots j_{k-1}j}^{i_{1}\cdots i_{k-1}i}h_{i_{1}}^{j_{1}}\cdots h_{i_{k-1}}^{j_{k-1}}. \]

The following properties will be used in later calculations (see \cite{Reilly,Barbosa-Colares} etc.).
\begin{proposition}\label{prop:sigmak}
	The following identities hold:
	\begin{itemize}
		\item[(i)] $(\sigma_{k}(A))_{j}^{i}h_{i}^{j}=k\sigma_{k}(A)$,
		\item[(ii)] $(\sigma_{k}(A))_{j}^{i}\delta_{i}^{j}=(n-k+1)\sigma_{k-1}(A)$,
		\item[(iii)] $(\sigma_{k}(A))_{j}^{i}h_{i}^{m}h_{m}^{j}=\sigma_{1}(A)\sigma_{k}(A)-(k+1)\sigma_{k+1}(A)$,
		\item[(iv)] $\frac{\partial}{\partial x_{i}}(\sigma_{k}(A))_{j}^{i}=0$.
	\end{itemize}
\end{proposition}

\subsection{The $k$-convexity of the hypersurface}
\label{sec:convexity}

The Gårding’s cone is defined by
\[ \Gamma_{k}:=\left\{ \lambda\in\mathbb{R}^{n}|\sigma_{i}(\lambda)>0 \text{ for } i=1,2,...,k\right\}. \] 
We say that $A\in\Gamma_{k}$ if its eigenvalues $\lambda(A)\in\Gamma_{k}$.

We recall Newton-MacLaurin inequalities
\begin{equation}\label{eq:NM1}
	H_{k-1}H_{k+1}\leq (H_{k})^2
\end{equation}
and
\begin{equation}\label{eq:NM2}
	H_{k}^{\frac{1}{k}}\leq H_{k-1}^{\frac{1}{k-1}}\leq \cdots \leq H_{1},
\end{equation}
where the equality occurs if and only if $\lambda_{1}=\lambda_{2}=\cdots=\lambda_{n}$ for $\lambda\in\Gamma_{k}$.

We also know the operator 
\[ (\sigma_{k}(A))_{j}^{i}\nabla_{i}\nabla^{j}:=(\sigma_{k}(A))_{j}^{i}g^{jl}\nabla_{i}\nabla_{l} \]
is elliptic if $A\in\Gamma_{k}$ (see \cite{Barbosa-Colares} etc.).

We say that a hypersurface $M$ is \emph{$k$-convex}, if its shape operator $A\in\Gamma_{k}$ everywhere in $M$. It is guaranteed in our settings by the following proposition, which can be found in \cite{Alias-Malacarne}.
\begin{proposition}\label{prop:k-convex}
	Suppose $M$ is a connected, spacelike hypersurface in $\mathbb{R}^{n,1}$ with boundary $\partial M$ on the hyperplane $\mathbb{R}^n\times \{ c \}$. If $M$ is not entirely contained in the hyperplane $\mathbb{R}^n\times \{ c \}$, and the $k$-th mean curvature $H_{k}$ of $M$ is constant, then the shape operator $A\in\Gamma_{k}$ everywhere in $M$.
\end{proposition}

\begin{proof}
	From Remark \ref{rem:u<0} and $u=c$ on $\partial M$, we know function $\langle X,E_{n+1} \rangle=-u$ attains its maximum at some interior point of $M$, denoted by $p_{0}$. Then \[ 0\geq \nabla_{i}\nabla_{j}\langle X,E_{n+1} \rangle(p_{0})=h_{ij}(p_{0})\langle N,E_{n+1} \rangle(p_{0}), \] here the inequality means the matrix $(\nabla_{i}\nabla_{j}\langle X,E_{n+1} \rangle)$ is non-positive definite. Since $\langle N,E_{n+1} \rangle<0$, we know $ h_{ij}(p_{0})\geq 0$. This implies $H_{k}>0$ and the shape operator $A(p_{0})\in \Gamma_{k}$. 
	
	At the same time, the set $\left\{ p\in M|A(p)\in\Gamma_{k}\right\} $ is open and closed, due to the smoothness of $u$, Newton-MacLaurin inequality \eqref{eq:NM2} and $H_{k}>0$. Then $A\in\Gamma_{k}$ everywhere in $M$, since $M$ is connected.
\end{proof}

\section{Proof of Theorem \ref{thm:main}}
\label{sec:pf}

Now we discuss under the assumption of Theorem \ref{thm:main} in the case $\Pi=\mathbb{R}^{n}\times \{ c \}$. Without loss of generality, we can suppose $H_{k}=1$ via a scaling, if $M$ is not entirely contained in $P$. The boundary conditions imply 
\begin{equation}\label{eq:bdry}
	u=c \text{ and } \langle N,E_{n+1}\rangle=\theta_{0} \quad\text{on }\partial\Omega,
\end{equation}
where $c$ and $\theta_{0}$ are constants. 

From now, let $\sigma_{k}=\sigma_{k}(A)$ and $(\sigma_{k})_{j}^{i}=(\sigma_{k}(A))_{j}^{i}$ for convenience. We establish an integral equality at first.
\begin{lemma}\label{lem:int}
	The following equality holds
	\begin{equation*}
		k\binom{n}{k}\int_{\Omega}(u-c)dx+(n-k+1)\int_{\Omega}\left(\langle N,E_{n+1}\rangle-\theta_{0} \right)\sigma_{k-1}dx=0.\label{eq:inteq}
	\end{equation*}
\end{lemma}

\begin{proof}
	Using Proposition \ref{prop:sigmak} (ii) and (iv), we have
	\begin{align*}
		&(n-k+1)\int_{\Omega}\sigma_{k-1}(\langle N,E_{n+1}\rangle-\theta_{0}) dx =\int_{\Omega}(\sigma_{k})_{j}^{i}\frac{\partial x_{j}}{\partial x_{i}}(\langle N,E_{n+1}\rangle-\theta_{0}) dx\\
		&\qquad =\int_{\Omega}\frac{\partial}{\partial x_{i}}\left((\sigma_{k})_{j}^{i}x_{j}(\langle N,E_{n+1}\rangle-\theta_{0})\right)dx-\int_{\Omega}(\sigma_{k})_{j}^{i}x_{j}\frac{\partial}{\partial x_{i}}\langle N,E_{n+1}\rangle dx.
	\end{align*}
	
	Then the boundary condition $\langle N,E_{n+1}\rangle=\theta_{0}$ and the divergence theorem imply
	\begin{equation}\label{eq:intO}
		(n-k+1)\int_{\Omega}\sigma_{k-1}(\langle N,E_{n+1}\rangle-\theta_{0}) dx=-\int_{\Omega}(\sigma_{k})_{j}^{i}x_{j}\frac{\partial}{\partial x_{i}}\langle N,E_{n+1}\rangle dx.
	\end{equation}
	
Using \eqref{eq:partialiN}, $(\sigma_{k})_{j}^{i}h_{i}^{l}=(\sigma_{k})_{i}^{l}h_{j}^{i}$ and the divergence theorem, we have
	\begin{align*}
		&-\int_{\Omega}(\sigma_{k})_{j}^{i}x_{j}\frac{\partial}{\partial x_{i}}\langle N,E_{n+1}\rangle dx  =\int_{\Omega}(\sigma_{k})_{j}^{i}x_{j}h_{i}^{l}u_{l}dx\\
		&\qquad =\int_{\Omega}(\sigma_{k})_{i}^{l}x_{j}h_{j}^{i}\frac{\partial}{\partial x_{l}}(u-c)dx =-\int_{\Omega}\frac{\partial}{\partial x_{l}}\left((\sigma_{k})_{i}^{l}x_{j}h_{j}^{i}\right)(u-c)dx.
	\end{align*}
	
	Since
	\begin{align*}
		&\frac{\partial}{\partial x_{l}}\left((\sigma_{k})_{i}^{l}x_{j}h_{j}^{i}\right)=(\sigma_{k})_{i}^{l}\frac{\partial}{\partial x_{l}}\left(x_{j}h_{j}^{i}\right)\\
		&\qquad =(\sigma_{k})_{i}^{j}h_{j}^{i}+x_{j}(\sigma_{k})_{i}^{l}\frac{\partial}{\partial x_{j}}h_{l}^{i}=k\sigma_{k}+x_{j}\frac{\partial}{\partial x_{j}}\sigma_{k}
	\end{align*}
	and $\sigma_{k}=\binom{n}{k}$, we know
	\begin{equation}\label{eq:-intO}
		-\int_{\Omega}(\sigma_{k})_{j}^{i}x_{j}\frac{\partial}{\partial x_{i}}\langle N,E_{n+1}\rangle dx=-k\binom{n}{k}\int_{\Omega}(u-c)dx.
	\end{equation}
	
	Combining \eqref{eq:intO} and \eqref{eq:-intO}, we finish the proof.
\end{proof}

Now we proof the theorem. We consider an auxiliary function \[ P=\langle X,E_{n+1} \rangle-\langle N,E_{n+1}\rangle. \]
By direct computation, we have
\[
\nabla_{i}P=\frac{\partial}{\partial x_{i}}P=\langle \frac{\partial X}{\partial x_{i}},E_{n+1} \rangle+h_{i}^{l}u_{l}
\]
and
\begin{align*}
	\frac{\partial^{2}P}{\partial x_{i}\partial x_{j}} & =\langle \frac{\partial^2 X}{\partial x_{i}\partial x_{j}},E_{n+1} \rangle+\frac{\partial}{\partial x_{j}}h_{i}^{l}u_{l}+h_{i}^{l}u_{lj}.
\end{align*}
Moreover,
\begin{align*}
	\nabla_{i}\nabla_{j}P & =\frac{\partial^{2}P}{\partial x_{i}\partial x_{j}}-\Gamma_{ij}^{m}\frac{\partial P}{\partial x_{m}}\\
	& =\langle \frac{\partial^2 X}{\partial x_{i}\partial x_{j}}-\Gamma_{ij}^{m}\frac{\partial X}{\partial x_{m}},E_{n+1} \rangle+\frac{\partial}{\partial x_{j}}h_{i}^{l}u_{l}+h_{i}^{l}u_{lj}-\Gamma_{ij}^{m}h_{m}^{l}u_{l}\\
	& =\langle\nabla_{i}\nabla_{j}X,E_{n+1}\rangle+\nabla_{j}h_{i}^{l}u_{l}-\Gamma_{mj}^{l}h_{i}^{m}u_{l}+h_{i}^{l}u_{lj}\\
	& =h_{ij}\langle N,E_{n+1}\rangle+\nabla_{j}h_{i}^{l}u_{l}-h_{i}^{m}\langle\nabla_{m}\nabla_{j}X,E_{n+1}\rangle\\
	& =(h_{ij}-h_{i}^{m}h_{mj})\langle N,E_{n+1}\rangle+\nabla_{j}h_{i}^{l}u_{l}.
\end{align*}

Then we obtain
\begin{align*}
	(\sigma_{k})_{j}^{i}\nabla_{i}\nabla^{j}P & =(\sigma_{k})_{j}^{i}g^{lj}\nabla_{i}\nabla_{l}P\\
	& =(\sigma_{k})_{j}^{i}(h_{i}^{j}-h_{i}^{m}h_{m}^{j})\langle N,E_{n+1}\rangle+(\sigma_{k})_{j}^{i}\nabla^{j}h_{i}^{l}u_{l}\\
	& =(k\sigma_{k}-\sigma_{1}\sigma_{k}+(k+1)\sigma_{k+1})\langle N,E_{n+1}\rangle+\nabla^{l}\sigma_{k}u_{l}.
\end{align*}

The Newton-MacLaurin inequalities and $H_{k}=1$ imply $\nabla^{l}\sigma_{k}=0$ and 
\[ \sigma_{1}\sigma_{k}-(k+1)\sigma_{k+1}\geq\frac{k}{n}\sigma_{1}\sigma_{k}\geq kH_{k}^{\frac{1}{k}}\sigma_{k}=k\sigma_{k}. \]
Here the equality occurs if and only if $\lambda_{1}=\cdots=\lambda_{n}=1$.

Combining these with $\langle N,E_{n+1}\rangle=-\frac{1}{\sqrt{1-|Du|^{2}}}<0$, we have \[ (\sigma_{k})_{j}^{i}\nabla_{i}\nabla^{j}P\geq0. \] Proposition \ref{prop:k-convex} ensures the operator $(\sigma_{k})_{j}^{i}\nabla_{i}\nabla^{j}$ is elliptic. Then the strong maximum principle gives either $P<-c-\theta_{0}$ in $\Omega$, or $P\equiv -c-\theta_{0}$ in $\overline{\Omega}$. 

In fact, the former can be excluded. Lemma \ref{lem:int} implies
\begin{align*}
	\int_{\Omega}(P+c+\theta_{0})\sigma_{k-1}dx&=-\int_{\Omega}(u-c)\sigma_{k-1}dx-\int_{\Omega}(\langle N,E_{n+1}\rangle-\theta_{0})\sigma_{k-1}dx\\
	&=-\int_{\Omega}(u-c)\left(\sigma_{k-1}-\binom{n}{k-1}\right)dx.
\end{align*}
We notice
\begin{align*}
	(\sigma_{k})_{j}^{i}\nabla_{i}\nabla^{j}u&=-(\sigma_{k})_{j}^{i}\langle \nabla_{i}\nabla^{j}X,E_{n+1} \rangle \\
	&=-(\sigma_{k})_{j}^{i}h_{i}^{j}\langle N,E_{n+1}\rangle=-k\sigma_{k}\langle N,E_{n+1}\rangle>0.
\end{align*}
The maximum principle implies $u<c$ in $\Omega$. The Newton-MacLaurin
inequality shows
\[ \sigma_{k-1}\geq\binom{n}{k-1}H_{k}^{\frac{k-1}{k}}=\binom{n}{k-1}.\]
Consequently,
\[ \int_{\Omega}(P+c+\theta_{0})\sigma_{k-1}dx \geq 0, \]
which contradicts $P<-c-\theta_{0}$ in $\Omega$.

Hence $P$ is constant in $\overline{\Omega}$, which implies 
\[\frac{\partial}{\partial x_{i}}\left(\frac{u_{j}}{\sqrt{1-|Du|^{2}}}\right)=h_{i}^{j}=\delta_{ij}.\]
As a result, $M$ is a part of a hyperboloid. Precisely,
\[ u=c+\theta_{0}+\sqrt{1+|x-a|^{2}} \] and $\Omega=B_{R}(a)$ is a ball,
where $R=\sqrt{\theta_{0}^{2}-1}$ and fixed $a\in\mathbb{R}^{n}$.


\begin{thebibliography}{00}

\bibitem{Aleksandrov}
Aleksandr Danilovich Aleksandrov, 
Uniqueness theorems for surfaces in the large. V, 
Vestnik Leningrad. Univ. {\bf 13} (1958), no.~19, 5--8.

\bibitem{Alias-Lopez-Pastor}
Luis~J. Al\'ias, Rafael L\'opez and Jos\'e A. Pastor, 
Compact spacelike surfaces with constant mean curvature in the Lorentz-Minkowski $3$-space, 
Tohoku Math. J. (2) {\bf 50} (1998), no.~4, 491--501.

\bibitem{Alias-Malacarne}
Luis~J. Al\'ias and J. Miguel Malacarne, 
Spacelike hypersurfaces with constant higher order mean curvature in Minkowski space-time, 
J. Geom. Phys. {\bf 41} (2002), no.~4, 359--375.

\bibitem{Alias-Pastor98}
Luis~J. Al\'ias and Jos\'e A. Pastor,
Constant mean curvature spacelike hypersurfaces with spherical boundary in the Lorentz-Minkowski space, 
J. Geom. Phys. {\bf 28} (1998), no.~1-2, 85--93.

\bibitem{Barbosa-Colares}
Jo\~ao Lucas Marques Barbosa and Ant\^onio Gervasio Colares, 
Stability of hypersurfaces with constant $r$-mean curvature, 
Ann. Global Anal. Geom. {\bf 15} (1997), no.~3, 277--297.

\bibitem{BNST}
Barbara Brandolini, Carlo Nitsch, Paolo Salani and Cristina Trombetti, 
Serrin-type overdetermined problems: an alternative proof, 
Arch. Ration. Mech. Anal. {\bf 190} (2008), no.~2, 267--280.

\bibitem{Cheng-Yau}
Shiu Yuen Cheng and Shing Tung Yau, 
Maximal space-like hypersurfaces in the Lorentz-Minkowski spaces, 
Ann. of Math. (2) {\bf 104} (1976), no.~3, 407--419.

\bibitem{GMY}
Shanze Gao, Hui Ma and Mingxuan Yang, 
Overdetermined problems for fully nonlinear equations with constant Dirichlet boundary conditions in space forms, 
Calc. Var. Partial Differential Equations {\bf 62} (2023), no.~6, Paper No. 183, 19 pp.

\bibitem{Hopf}
Heinz Hopf, 
Differential geometry in the large, second edition, 
Lecture Notes in Mathematics, 1000, Springer, Berlin, 1989.

\bibitem{Jia}
Xiaohan Jia,
Overdetermined problems for Weingarten hypersurfaces, 
Calc. Var. Partial Differential Equations {\bf 59} (2020), no.~2, Paper No. 78, 15 pp.

\bibitem{LiAnMin}
An Min Li, 
Spacelike hypersurfaces with constant Gauss-Kronecker curvature in the Minkowski space, 
Arch. Math. (Basel) {\bf 64} (1995), no.~6, 534--551.

\bibitem{Lopez}
Rafael L\'opez,
Constant mean curvature surfaces with boundary,
Springer Monographs in Mathematics, Springer, Heidelberg, 2013.

\bibitem{Montiel-Ros}
Sebastián Montiel and Antonio Ros, 
Compact hypersurfaces: the Alexandrov theorem for higher order mean curvatures, 
{\it Differential geometry}, 279--296, Pitman Monogr. Surveys Pure Appl. Math., 52, Longman Sci. Tech., Harlow.

\bibitem{Reilly}
Robert~C. Reilly, 
On the Hessian of a function and the curvatures of its graph, 
Michigan Math. J. {\bf 20} (1973), 373--383.

\bibitem{Serrin}
James Serrin, 
A symmetry problem in potential theory, 
Arch. Rational Mech. Anal. {\bf 43} (1971), 304--318.

\bibitem{Treibergs}
Andrejs~E. Treibergs, 
Entire spacelike hypersurfaces of constant mean curvature in Minkowski space, 
Invent. Math. {\bf 66} (1982), no.~1, 39--56.

\bibitem{Wang-Xiao}
Zhizhang Wang and Ling Xiao, 
Entire spacelike hypersurfaces with constant $\sigma_k$ curvature in Minkowski space, 
Math. Ann. {\bf 382} (2022), no.~3-4, 1279--1322.

\bibitem{Weinberger}
Hans F. Weinberger,
Remark on the preceding paper of the Serrin,
Arch. Rational Mech. Anal. {\bf 43} (1971), 319--320.

\bibitem{Wente}
Henry~C. Wente, 
The symmetry of sessile and pendent drops, 
Pacific J. Math. {\bf 88} (1980), no.~2, 387--397.

\end{thebibliography}
\end{document}